\documentclass[12pt,legno]{article}
\usepackage{verbatim, amsmath,amssymb,amsthm,amscd,euscript,mathrsfs}
\usepackage{fancyhdr}
\input xy
\xyoption{all}

\newtheorem{theorem}{Theorem}

\newtheorem{proposition}[theorem]{Proposition}

\newtheorem{corollary}[theorem]{Corollary}

\theoremstyle{definition}
\newtheorem{definition}[theorem]{Definition}

\def\hom{\mathop\mathrm{Hom}\nolimits}

\def\rank{\mathop\mathrm{rank}\nolimits}

 \def\mP{\mathbb P}

 \def\kN{\mathcal N}
 \def\kO{\mathcal O}
 \def\kP{\mathcal P}
 
\def\kE{\mathcal E} 
\def\kF{\mathcal F}

\def\kI{\mathcal I}

\begin{document}

\title{On irreducibility of the family of ACM curves of degree 8 and genus 4 in $\mP^4_k$}
\author{Elena Drozd}
\date{}

\maketitle

\begin{abstract}
Let $C$ be an arithmetically Cohen-Macaulay curve of arithmetic genus 
4. We prove that the family of such curves of degree 8 in $\mP^4_k$ is irreducible.
\end{abstract}


~\newline

Liaison using complete intersection or Gorenstein schemes is widely used in 
algebraic geometry. An excellent reference book is \cite{migliore}. We use the technique
of resolving ideal sheaves of ACM curves
by special type of sheaves (so called $\kE$- and $\kN$-type resolutions) to link
the curve in question to a simpler curve.  
Using this technique we conclude that a family of $(8,4)$ curves is 
irreducible.  Also this paper demonstrates a usage of correspondance between ACM
curves and ACM sheaves ( as discussed in \cite{drozd1})

For convenience we recall here some definitions and results of liaison
theory we will be using in this work. See \cite{migliore} for  reference. 
\begin{definition}
A scheme $X$ of $\mP^n_k$ is called arithmetically Cohen-Macaulay (ACM)
if its homogeneous coordinate ring is a Cohen-Macaulay ring.
\end{definition}
\begin{definition}
Let $Z$ be a subscheme of $\mP^4_k$. Let $V_1, V_2$ be equidimensional 
subschemes of $\mP^n_k$ of codimension $r$ and without embedded
components. We say that $V_1$ and $V_2$ are linked by $Z$ if
\begin{enumerate}
\item $\kI_Z\subset \kI_{V_1}\cap \kI_{V_2}$
\item $\kI_{V_2}/\kI_Z \cong \hom_{\kO_{\mP^n_k}}(\kO_{V_1},\kO_Z)$
\item $\kI_{V_1}/\kI_Z \cong \hom_{\kO_{\mP^n_k}}(\kO_{V_2},\kO_Z)$
\end{enumerate}
If $Z$ is AG, we say $V_1$ is G-linked to $V_2$;
if $Z$ is CI, we say $V_1$ is CI-linked to $V_2$.
\end{definition}

\begin{definition}
On a nonsingular quadric hypersurface $Q$ a locally free sheaf $\kF$ with the property that $H^i_*(\kF)=0$ for $i=1,2$ is 
called and \emph{ACM sheaf}.
\end{definition}

\begin{proposition}\label{p81}
Let $C$ be an ACM curve of degree 8 and arithmetic genus 4 in $\mP^4_k$. Then $\kI_C$ is generated
in degrees 2 and 3.
\end{proposition}
\begin{proof}
We need to compute the cohomology table of $\kI_C(n)$. From the
 Rieman-Roch theorem 
 $h^0\bigl(\kO_C(1)\bigr)-h^1\bigl(\kO_C(1)\bigr)=8+1-4=5$.  
Taking chomology in the exact sequence $$0\longrightarrow \kI_C(n)\longrightarrow \kO_{\mP^4_k}(n) \longrightarrow \kO_C(n)\longrightarrow 0$$ we arrive
at 
\begin{multline}\label{cohom1}
0\longrightarrow 
H^0\bigl(\kI_C(n)\bigr) \longrightarrow 
H^0\bigl(\kO_{\mP^4_k}(n)\bigr) \longrightarrow 
H^0\bigl(\kO_C(n)\bigr) \longrightarrow \\ \longrightarrow 
H^1\bigl(\kI_C(n)\bigr) \longrightarrow 
H^1\bigl(\kO_{\mP^4_k}(n)\bigr) \longrightarrow 
H^1\bigl(\kO_C(n)\bigr) \longrightarrow \\ \longrightarrow
H^2\bigl(\kI_C(n)\bigr) \longrightarrow 
H^2\bigl(\kO_{\mP^4_k}(n)\bigr) \longrightarrow 
H^2\bigl(\kO_C(n)\bigr) \longrightarrow 
H^3\bigl(\kI_C(n)\bigr) \longrightarrow \dots \quad.
\end{multline}
$H^1\bigl(\kI_C(n)\bigr)=0$ since $C$ is ACM. Thus, map $H^0\bigl(\kO_{\mP^4_k}(n)\bigr)\longrightarrow H^0\bigl(\kO_C(n)\bigr)$ is 
surjective and so, $h^0\bigl(\kO_C(n)\bigr)\leq h^0\bigl(\kO_{\mP^4_k}(n)\bigr)$. From the Rieman-Roch theorem
we get 

\noindent $h^0\bigl(\kO_C(1)\bigr)=5+h^1\bigl(\kO_C(1)\bigr)$ while 
$h^0\bigl(\kO_{\mP^4_k}(1)\bigr)=5$. 
This implies
that $h^1\bigl(\kO_C(1)\bigr)=0$ and $h^0\big(\kO_C(1)\big)=5$.

From the exact sequence (\ref{cohom1}) we obtain $h^2\bigl(\kI_C(n)\bigr)=h^1\bigl(\kO_C(n)\bigr)$
since $H^i\bigl(\kO_{\mP^4_k}(n)\bigr)=0$ for $i=1,2$ \cite[III.5.1]{hartshorne}. Wherefrom $h^2\bigl(\kI_C(1)\bigr)=h^1\bigl(\kO_C(1)\bigr)=0$. 
Note also that $H^3(\kI_C)=0$. 
Thus, the cohomology table is:
$$
   
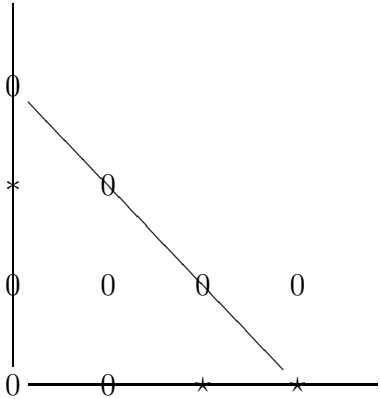
\xymatrix{ 
 &  {}\ar@{-}[dddd] \\
 & 0 \ar@{-}[dddrrr]\\
 & {*} & 0 \\
 & 0&0&0&0 \\
 & 0\ar@{-}[rrrr] &0&\star &\star & & & . }
 $$ 
 Thus, $h^i\bigl(\kI_C(3-i)\bigr)=0$ for all $i>0$. By definition 1.1.4 of \cite{migliore} this implies
 that $\kI_C$ is 3-regular. This, in turn, by Castelnuovo-Mumford regularity \cite{migliore} 1.1.5.(1), implies
 that $h^i\bigl(\kI_C(k)\bigr)=0$ for $i>0$, $k+i\geq 3$. Equivalently, 
 $h^2\bigl(\kI_C(n)\bigr)=0\quad \text{ for }n\geq 1.$
 Thus, $h^1\bigl(\kO_C(n)\bigr)=0$ for all $n\geq 1$.
 Therefore, $\kI_C(k)$ is generated as $\kO_{\mP^4_k}$-module by its global sections for all $k\geq 3$ 
 (by \cite{migliore} theorem 1.1.5.(3) ).
 \end{proof}

\begin{corollary}
Any ACM curve $C$ of degree 8 and genus 4 in $\mP^4_k$ is contained in a quadric hypersurface.
\end{corollary}
\begin{proof}
 
 Proposition above implies that
 $h^0(\kO_C(n))=nd+1-g$
 for $n\ge1$, or

 $$ 
    \begin{array}{c|c|c|c|c|c}
  n&0&1&2&3&4\\
 \hline
 h^0(\kO_C(n))&1&5&13&21&29 .
 \end{array} 
 $$ 
 Recall that we have 
 $$ 
    \begin{array}{c|c|c|c|c|c}
  n&0&1&2&3&4\\
\hline
 h^0(\kO_Q(n))&1&5&15&35&70,
 \end{array} 
 $$ 
wherefrom we obtain 
 $$ 
    \begin{array}{c|c|c|c|c|c}
  n&0&1&2&3&4\\
\hline
 h^0(\kI_C(n))&0&0&2 &14 &41 .
 \end{array} 
 $$ 
Thus $h^0\big(\kI_C(2)\big)=2$, which implies that there is at least one quadric hypersurface
containing $C$.
\end{proof}

\begin{proposition}
There is no ACM curve of degree 8 and genus 4 in $\mP^3_k$. 
\end{proposition}
\begin{proof}
Assume $C$ is an ACM curve of degree 8 and genus 4 in $\mP^3_k$. Taking cohomology in the 
exact sequence $$0\longrightarrow \kI_C(1)\longrightarrow \kO_{\mP^3_k}(1)\longrightarrow \kO_C(1) \longrightarrow 0$$ 
we obtain:
\begin{equation*}
0\longrightarrow H^0\big(\kI_C(1)\big)\longrightarrow H^0\big(\kO_{\mP^3_k}(1)\big) \longrightarrow H^0\big(\kO_C(1)\big) \longrightarrow H^1\big(\kI_C(1)\big) \longrightarrow 0 .
\end{equation*}
Thus $h^0\big(\kI_C(1)\big)=h^0\big(\kO_{\mP^3_k}(1)\big)-h^0\big(\kO_C(1)\big)$ since $h^1\big(\kI_C(1)\big)=0$
for an ACM curve $C$. However $h^0\big(\kO_{\mP^3_k}(1)\big)=4$ and $h^0\big(\kO_C(1)\big)=5+h^1\big(\kO_C(1)\big)\geq 5$.
This would give $h^0\big(\kI_C(1)\big)<0$, which is impossible. Thus $h^1\big(\kI_C(1)\big)\neq 0$, and 
$C$ is not an ACM curve.
\end{proof}

\begin{corollary}\label{ndg84}
Any ACM curve of degree 8 and genus 4 in $\mP^4_k$ is nondegenerate.
\end{corollary}

\begin{proposition}\label{p94}
Let $C$ be an ACM curve of degree 8 and genus 4 on a nonsingular quadric hypersurface $Q$. Then, 
there is an $\kE$-type resolution  of $\kI_C$ of the form 
\begin{equation}\label{equation:e91}
0\longrightarrow \kE_0^2(-2)\longrightarrow \kO_Q(-2)\oplus\kO_Q^4(-3)\longrightarrow\kI_C\longrightarrow0
\end{equation}
\begin{proof}
Since $0\longrightarrow \kI_C(n)\longrightarrow \kO_Q(n)\longrightarrow \kO_C(n)\longrightarrow 0$ is exact we obtain
 $$ 
    \begin{array}{c|c|c|c|c|c|c|c}
  n&0&1&2&3&4&5&6\\
\hline
 h^0\big(\kI_C(n)\big)&0&0&1 &9 &26 &54 & 95 .
 \end{array} 
 $$ 

We know that $\kI_C$ is generated in degrees 2 and 3 and also 
the generator of $\kI_C$ in degree 2 multiplied by linear functions
 gives a 5-dimensional subspace of $H^0\bigl(Q,\kI_C(3)\bigr)$. Therefore we need 4 generators
in degree 3, which are not products of linear form and the degree two generator. 
Thus, 
there is  an $\kE$-type resolution of $\kI_C$ of the form
  \begin{equation}\label{equation:eType84}
   0\longrightarrow \kE\longrightarrow \kO_Q(-2)\oplus \kO_Q^4(-3)\longrightarrow \kI_C\longrightarrow 0
  \end{equation}
where $\kE$ is ACM sheaf by \cite[Theorem 2]{drozd1}  and $\rank \kE=4$.

(\ref{equation:eType84}) gives the following table of cohomology:
 $$ 
    \begin{array}{c|c|c|c}
  n&h^0(\kE(n))& h^0\bigl(\kO_Q(-2+n)\oplus\kO_Q^4(-3+n)\bigr)&h^0(\kI_C(n)\\
 \hline
 0&0&0&0\\
\hline
1&0&0&0\\
\hline
2&0&1&1\\
\hline
3&0&9&9\\
\hline
4&8&34&26\\
\hline
5&32&86&54\\
\hline
6&80&175&95 \\
\hline
 \end{array} 
 $$ 

Thus, by  \cite[Corollary 3]{drozd1} $\kE$ must be one of the following:
\begin{enumerate}
\item $\kE_0(a)\oplus \kO(b) \oplus \kO(c)$, or
\item $\kE_0(a)\oplus \kE_0(b)$, or
\item $\bigoplus_{i=1}^4\kO(a_i)$,
\end{enumerate}
where the sheaf $\kE_0$ is given by cite[Definition 5]{drozd1} . 

Comparing 
 $$ 
  \begin{array}{c|c|c|c|c|c|c|c} 
   n&0&1&2&3&4&5&6\\
\hline
h^0(\kE(n))&0&0&0&0&8&32&80\\
\end{array} 
 $$
 with cohomology tables for $\kE_0(n)$ and $\kO(n)$ we obtain 
 $\kE=\kE_0^2(-2)$. 
This gives us $$0\longrightarrow \kE_0^2(-2)\longrightarrow \kO_Q(-2)\oplus\kO_Q^4(-3)\longrightarrow\kI_C\longrightarrow 0 $$
as an $\kE$-type resolution of $\kI_C$.
 \end{proof}
 \end{proposition} 

%

  \begin{proposition}\label{p93}
 Let $C$ be an ACM curve of degree 8 and genus 4 on $Q$. Then $C$ 
 is CI-linked to an ACM curve $C^\prime$ of degree
 4 and genus 0 (possibly reducible).  
 
 \begin{proof}
Let $C$ be an ACM curve of degree 8 and genus 4 on a quadric hypersurface $Q$ in $\mP^4_k$.

Then by \ref{ndg84} $C$ is nondegenerate. Note that $h^0\bigl(P^4_k,\kI_C(2)\bigr)=2$, therefore
$h^0\bigl(Q,\kI_C(2)\bigr)=1$. Thus 
the generator of $\kI_C$ in degree 2 cuts out a surface $Y$ of degree 4 on $Q$. We claim that $Y$ is 
irreducible. To prove this let $Y$ be a union of two  surfaces $Q_1$ and $Q_2$. Then $\deg Q_1 =\deg Q_2=2$ since
$Y$ is a degree 4 surface on a nonsingular quadric hypersurface and thus Klein's theorem \cite[ex.II.6.5.(d)]{hartshorne}
implies that $\deg Q_i,\quad i=1,2$ must be even. However a quadric surface lies in $\mP^3_k$, which contradicts
$C$ is nondegenerate. Thus $Y$ is irreducible.

Let $F$ be a hypersurface of degree 3 containing $C$, but not containing $Y$ completely. Such a hypersurface
exists since $h^0\big(Q,\kI_C(3)\big)-\dim V=14-5=9$, where $V$ is a subspace of $H^0\big(Q,\kI_C(3)\big)$
generated by elements of the form $l\cdot s$, where $l$ is a linear form and $s\in H^0\big(Q,\kI_C(2)\big)$.
Let $Z$ be a complete intersection of $Y$ and $F$. Then $Z$ has degree 12 and it contains $C$. Let curve
$C'$ be CI-linked to $C$ via $Z$. 

%
Note that $C^\prime$ is ACM since so is $C$. To complete the proof we need to
compute degree and genus of $C^\prime$: 
$\deg C^\prime = \deg Z-\deg C=4$. By \cite[corollary 5.2.14]{migliore},
$g(C)-g(C^\prime)=\frac12(\deg F+\deg Y-5)\cdot (\deg C-\deg C^\prime)$. Thus, $g(C^\prime)=g(C)-4=0$,
wherefrom $C^\prime$ is an ACM curve of degree 4 and genus 0. 
 \end{proof}
 \end{proposition} 
%
In order to find an $\kN$-type resolution of a nondegenerate ACM (8,4) curve, we will
determine an $\kE$-type resolution
 of a linked (4,0) ACM curve.

Now we compute an $\kE$-type resolution of a nondegenerate ACM (4,0) curve $C^\prime$. 

 \begin{proposition}\label{etexists}
There exists an $\kE$-type resolution of any 
ACM (4,0) curve on a nonsingular quadric hypersurface $Q$  of the form 
 $$ 
   0\longrightarrow\kE^2_0(-1)\longrightarrow\kO^5_Q(-2)\longrightarrow\kI_C\longrightarrow0.
 $$
 \begin{proof}
We claim that $\kI_C$ is generated in degree 2 and $h^1\bigl(\kO_C(n)\bigr)=0$ for $n\geq 1$.

To prove this we need to compute the cohomology table. From the Rieman-Roch theorem 
$h^0\bigl(\kO_C(1)\bigr)-h^1\bigl(\kO_C(1)\bigr)=4+1-0=5$. Thus $h^0\bigl(\kO_C(1)\bigr)\geq 5$. 
Taking cohomology in the short exact sequence $$0\longrightarrow \kI_C(n)\longrightarrow \kO_Q(n) \longrightarrow \kO_C(n)\longrightarrow 0$$ we 
obtain:
\begin{multline}\label{cohom2}
0\longrightarrow 
H^0\bigl(\kI_C(n)\bigr) \longrightarrow 
H^0\bigl(\kO_Q(n)\bigr) \longrightarrow 
H^0\bigl(\kO_C(n)\bigr) \longrightarrow \\ \longrightarrow 
H^1\bigl(\kI_C(n)\bigr) \longrightarrow 
H^1\bigl(\kO_Q(n)\bigr) \longrightarrow 
H^1\bigl(\kO_C(n)\bigr) \longrightarrow \\ \longrightarrow 
H^2\bigl(\kI_C(n)\bigr) \longrightarrow 
H^2\bigl(\kO_Q(n)\bigr) \longrightarrow 
H^2\bigl(\kO_C(n)\bigr) \longrightarrow 
H^3\bigl(\kI_C(n)\bigr) \longrightarrow \dots . 
\end{multline}

Note that $H^1\bigl(\kI_C(n)\bigr)=0$ since $C$ is ACM. Thus the map 
$H^0\bigl(\kO_Q(n)\bigr) \longrightarrow 
H^0\bigl(\kO_C(n)\bigr)$ is surjective, and $h^0\bigl(\kO_C(n)\bigr)\leq h^0\bigl(\kO_Q(n)\bigr)$.
For $n=0$ this means that $h^0(\kO_C)\leq h^0(\kO_Q)=1$. However, $h^0(\kO_C)-h^1(\kO_C)=1$, thus 
$h^0(\kO_C)=1$ and $h^1(\kO_C)=0$.

Also, since $h^1\bigl(\kO_Q(n)\bigr)= h^2\bigl(\kO_Q(n)\bigr)= 0$ \cite[Ex.III.5.5(c)]{hartshorne} we have 
$h^1\bigl((\kO_C(n)\bigr)=h^2\bigl((\kI_C(n)\bigr)$. 
Thus $h^2(\kI_C)=h^1(\kO_C)=0$. 
We obtain the following 
cohomology table:

$$
   
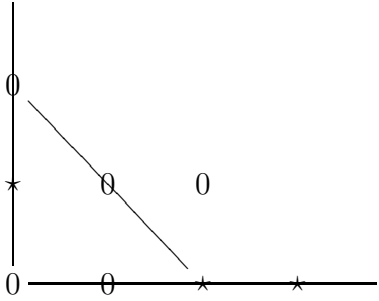
\xymatrix{ 
 &  {}\ar@{-}[ddd] \\
 & 0 \ar@{-}[ddrr]\\
 & \star & 0 & 0 &  \\
 & 0\ar@{-}[rrrr] &0&\star &\star & & & . }
 $$ 

Thus, by Castelnuovo-Mumford regularity $\kI_C$ is generated in degree 2 and $\kO_C(n)$ is nonspecial
for $n\geq 1$. 

Thus $ h^0\bigl((\kO_C(n)\bigr)=nd+1-g=4n+1$ for $n\geq 1$ or 
 \begin{align*}
  & \begin{array}{c|c|c|c|c|c|c|c}
  n&0&1&2&3&4&5&6\\
\hline
 h^0(\kO_C(n))&1&5&9&13&17&21&25
 \end{array}  \quad \text{ and } \\
    & \begin{array}{c|c|c|c|c|c|c|c}
  n&0&1&2&3&4&5&6\\
\hline
 h^0(\kO_Q(n))&1&5&14&30&55&91&140
 \end{array} \quad ,
 \end{align*} 
  which implies the following table 
 $$ 
      \begin{array}{c|c|c|c|c|c|c|c}
  n&0&1&2&3&4&5&6\\
\hline
 h^0(\kI_C(n))&0&0&5&17&38&70&115
 \end{array} 
 $$ 
 since $0\longrightarrow\kI_C(n)\longrightarrow\kO_Q(n)\longrightarrow\kO_C(n)\longrightarrow0$ is exact.

Thus we have the following exact sequence: 
 $$ 
   0\longrightarrow\kE\longrightarrow\kO^5_Q(-2)\longrightarrow\kI_C\longrightarrow0
 $$ 
 with $\kE$ ACM sheaf of rank 4. Thus, by \cite[Corollary 3]{drozd1} $\kE$ is one of the 
 following: 
 \begin{enumerate}
\item $\kE_0(a)\oplus \kO(b) \oplus \kO(c)$, or
\item $\kE_0(a)\oplus \kE_0(b)$, or
\item $\bigoplus_{i=1}^4\kO(a_i)$.
\end{enumerate}
$h^0\bigl(\kE(n)\bigr)=h^0\bigl(\kO_Q(n-2)\bigr)- h^0\bigl(\kI_C(n)\bigr)$ since $h^1(\kE)=0$. 
Comparing cohomology table: 
 $$ 
      \begin{array}{c|c|c|c|c|c|c|c}
  n&0&1&2&3&4&5&6\\
\hline
 h^0(\kE(n))&0&0&0&8 &32&80&160
 \end{array} 
 $$  with cohomology tables $\kE_0(n)$ and $\kO(n)$ we obtain
   $\kE=\kE_0^2(-1)$, proving the proposition. 
 \end{proof} 
 \end{proposition} 

This proposition together with \cite[Corollary 1]{drozd1}  give us the following

\begin{corollary}\label{c40}
All ACM (4,0) curves on a nonsingular quadric hypersurface $Q$ form an irreducible family.
\end{corollary}

 \begin{proposition}
  There exists an $\kN$-type resolution of an ACM (8,4) curve $C$ on a nonsingular quadric
  hypersurface $Q$ in $\mP^4_k$ of the form
 $$ 
   0\longrightarrow\kO_Q^5(-5)\longrightarrow\kO_Q(-4)\oplus\kO_Q(-3)\oplus \kE_0^2(-3)\longrightarrow \kI_{C}\longrightarrow0.
 $$ 
 \begin{proof}
By proposition \ref{p93} an ACM curve of degree 8 and genus 4 on a nonsingular
quadric hypersurface $Q$ in $\mP^4_k$ can be CI-linked to an ACM curve $C^\prime$ of degree 4 and genus 0
by a complete intersection curve $Z$ formed by two divisors $\kO_Q(4)$ and $\kO_Q(3)$. 
By proposition \ref{etexists} there
exists
an $\kE$-type resolution of $\kI_{C^\prime}$ of the form:
 $$ 
   0\longrightarrow\kE^2_0(-1)\longrightarrow\kO^5_Q(-2)\longrightarrow\kI_{C^\prime}\longrightarrow0.
 $$
However, by \cite[Proposition 2]{drozd1},
 $$ 
   \big(\kE_0^2(-1)\big)^\vee=\big(\kE_0^\vee(1)\big)^2=\kE_0^2(4).
 $$ 
Thus, there exists an $\kN$-type resoluion of $\kI_C$ of the form:
 $$ 
   0\longrightarrow\kO_Q^5(-5)\longrightarrow\kO_Q(-4)\oplus\kO_Q(-3)\oplus \kE_0^2(-3)\longrightarrow \kI_{C}\longrightarrow0.
 $$ 
 \end{proof} 
 \end{proposition} 

The above proposition together with \cite[Corollary1]{drozd} imply

\begin{corollary}
All ACM (8,4) curves  on a nonsingular quadric hypersurface $Q$ form an irreducible family.
\end{corollary}

We note here that any nonsingular curve of degree 8 and genus 4 on $Q$ is ACM.

\begin{proposition}
Let $C$ be a nonsingular curve of degree 8 and genus 4 on a nonsingular quadric hypersurface $Q$ in $\mP^4_k$.
Then $C$ is ACM. 
\end{proposition}
\begin{proof}
$h^1\big(\kO_C(n)\big)=0$ for all $n\geq 1$ since $2g-2=6\leq \deg C$. By the Rieman-Roch theorem
$h^0\big(\kO_C(1)\big)=8+1-4+h^1\big(\kO_C(1)\big)$. Thus $h^0\big(\kO_C(1)\big)=5$. Taking 
cohomology in the exact sequence $0\longrightarrow \kI_C(1)\longrightarrow \kO_Q(1)\longrightarrow \kO_C(1)\longrightarrow 0$ we obtain
$$
0\longrightarrow H^0\big(\kI_C(1)\big)\longrightarrow H^0\big(\kO_Q(1)\big)\longrightarrow H^0\big(\kO_C(1)\big)\longrightarrow H^1\big(\kI_C(1)\big)\longrightarrow 0
$$
If $H^1\big(\kI_C(1)\big)\neq 0$ then $H^0\big(\kI_C(1)\big)\neq 0$. Therefore there exists a hyperplane
$H$ such that $C\subset H\cap Q$ and $H\cap Q$ is a surface of degree 2 in $\kP^3_k$. Thus $H\cap Q$ is one
of the following:
\begin{itemize}
\item ~ two planes, or 
\item ~ double plane, or 
\item ~ quadric cone, or 
\item ~ nonsingular quadric surface.
\end{itemize}
Two planes and double plane are impossible since there is no (8,4) curve in $\mP^2_k$ (plane curve of degree
8 has genus 21). The set of possible pairs $(d,g)$ on a quadric cone is a subset of the set of possible
pairs $(d,g)$ on a nonsingular quadric surface. But there is no (8,4) curve on a quadric surface in $\mP^3_k$.
Thus, $H^1\big(\kI_C(1)\big)$ must be zero. Similarly $H^1\big(\kI_C(2)\big)=0$. From the exact sequence
$$
0\longrightarrow H^0\big(\kI_C(2)\big)\longrightarrow H^0\big(\kO_Q(2)\big)\longrightarrow H^0\big(\kO_C(2)\big)\longrightarrow H^2\big(\kI_C(2)\big)\longrightarrow 0
$$
we obtain  $h^0\big(\kO_C(2)\big)=13$ and $h^0\big(\kO_Q(2)\big)=14$. Thus $h^0\big(\kI_C(2)\big)\geq 1$.
We claim that $h^0\big(\kI_C(2)\big)=1$.
If $h^0\big(\kI_C(2)\big)\geq 2$ then $h^0\big(\mP^4_k,\kI_C(2)\big)\geq 3$. Thus $C$ must be contained 
in the intersection $Z$ of three quadric surfaces, wherefrom $Z$ must be one of the following:
\begin{enumerate}
\item A curve. Then it is of degree 8 and genus 5, or
\item A surface of degree $\leq 4$. 
\end{enumerate}
Neither of these is possible, therefore $h^1\big(\kI_C(2)\big)=0$. Note that for any curve $C$ 
$h^0\big(\kO_Q\big)\cong k \cong h^0\big(\kO_C\big)$. Therefore $h^1\big(\kI_C\big)=0$. 
Also, $h^2\big(\kI_C(1)\big)=h^1\big(\kO_C(1)\big)=0$ and $h^3\big(\kI_C\big)=0$.
Thus we have the 
following cohomology table for $\kI_C(n)$:

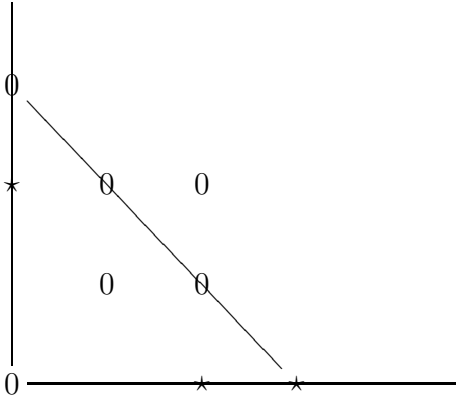
$$
   \xymatrix{ 
 &  {}\ar@{-}[dddd] \\
 & 0 \ar@{-}[dddrrr]\\
 & \star & 0 & 0 &  \\
 &  & 0 & 0 &  \\
 & 0\ar@{-}[rrrrr] & &\star &\star & & & . }
 $$ 
Thus $\kI_C$ is 2-regular and $h^1\big(\kI_C(n)\big)=0$ for $n\geq 2$ and for $n<0$. 
However $h^1\big(\kI_C\big)=h^1\big(\kI_C(1)\big)=h^1\big(\kI_C(2)\big)=0$, therefore $C$ is ACM.
\end{proof}

~\newline


\begin{thebibliography}{99}

\bibitem{cdh} M.~Casanellas, E.~Drozd, R.~Hartshorne: \emph{Gorenstein Liason and
 ACM Sheaves} accepted by Crelle's Journal, (2004)


\bibitem{drozd} E.~Drozd: \emph{Curves on a nonsingular quadric hypersurface 
in $\mP^4_k$: existence and liaison theory}, Ph.D thesis, UC Berkeley, (2003).

\bibitem{drozd1}E.~Drozd: \emph{ACM sheaves on a nonsingular Quadric hypersurface in 
$\mP^4_k$}, arXiv: math.AG/0409243 (2004)

\bibitem{eisenbud} D.~Eisenbud: \emph{Commuative Algebra with a View Toward Algebraic Geometry}, 
Springer (1999).

\bibitem{hartshorne} R.~Hartshorne: \emph{Algebraic Geometry}, Springer (1977).

\bibitem{knorrer} H.~Knorrer: \emph{Cohen-Macaulay modules on hypersurface singularities I}, Invent. Math. {\bf 88}
(1987) 153-164.


\bibitem{migliore} J.C.~Migliore: \emph{Gorenstein liaison theory and deficiency modules}, Progress in
Mathematics {\bf 165}  Birkh$\ddot{a}$user(1998).


\end{thebibliography}
\end{document}